\documentclass[review]{elsarticle}

\usepackage{amsmath,amssymb,bm}
\usepackage{physics}
\usepackage{color}
\usepackage{comment}
\usepackage{graphicx}

\usepackage{lineno,hyperref}
\usepackage{booktabs}
\usepackage{multirow}

\def\<{\langle}
\def\>{\rangle}

\newcommand{\ct}{\mathrm{\ast}} 

\newcommand{\R}{\mathop{\mathbb{R}}\nolimits}
\newcommand{\CA}{\mathop{\mathbb{C}}\nolimits}

\newcommand{\ov}[1]{\overline{#1}}

\renewcommand{\Re}{\mathop{\mathrm{Re}}}

\newtheorem{Thm}{Theorem}
\newtheorem{Prop}{Proposition}

\newtheorem{Lem}{Lemma}

\begin{document}

\begin{frontmatter}

\title{One parameter generalization of BW inequality and its application to open quantum dynamics}

\author[1]{Dariusz Chru\'sci\'nski}
\ead{darch@fizyka.umk.pl}
\author[2]{Gen Kimura}
\ead{gen@shibaura-it.ac.jp}
\author[3]{Hiromichi Ohno\corref{C}}
\ead{h\_ohno@shinshu-u.ac.jp}
\author[1]{Tanmay Singal}
\ead{tanmaysingal@gmail.com}

\address[1]{Institute of Physics, Faculty of Physics, Astronomy and Informatics Nicolaus Copernicus University, Grudzi\c{a}dzka 5/7, 87--100 Toru\'n, Poland}
\address[2]{College of Systems Engineering and Science,
Shibaura Institute of Technology, Saitama 330-8570, Japan}
\address[3]{Department of Mathematics, Faculty of Engineering, Shinshu University,
4-17-1 Wakasato, Nagano 380-8553, Japan.} 
\cortext[C]{Corresponding author}
\begin{abstract}
In this paper, we introduce a one parameter generalization of the famous B\"ottcher-Wenzel (BW) inequality in terms of a $q$-deformed commutator. 
For $n \times n$ matrices $A$ and $B$, we consider the inequality
\[
\Re\<[B,A],[B,A]_q\> \le c(q) \|A\|^2 \|B\|^2,
\]
where $\<A,B\> = \tr(A^*B)$ is the Hilbert-Schmidt inner product, $\|A\|$ is the Frobenius norm,
$[A,B] =AB-BA$ is the commutator, and $[A,B]_q =AB-qBA$ is the $q$-deformed commutator. 
We prove that when $n=2$, or when $A$ is normal with any size $n$, the optimal bound is given by
\[
c(q) = \frac{(1+q) +\sqrt{2(1+q^2)}}{2}.
\]
We conjecture that this is also true for any matrices, and this conjecture is perfectly supported for $n$ up to $15$ by numerical optimization.
When $q=1$, this inequality is exactly BW inequality.
When $q=0$, this inequality leads the sharp bound for the $r$-function which is recently derived for the application to universal constraints of relaxation rates in open quantum dynamics. 
\end{abstract}
\end{frontmatter}

\section{Introduction}

For any complex matrices $A,B \in M_n(\CA)$, the sharp bound for the norm of commutator is known: 
\begin{equation}\label{BWineq}
\|[A,B]\|^2 \le 2 \|A\|^2\|B\|^2,	
\end{equation}
where $[A,B] := AB-BA$ denotes the commutator and the matrix norm is the Frobenius norm $\|A\| = \sqrt{\tr (A^*A)}$. 
The inequality is often called the B\"ottcher-Wenzel (BW) inequality because of the following history: The inequality was first conjectured by B\"ottcher and Wenzel in \cite{ref:BW} giving a proof for real $2\times 2$ matrices and also for normal matrices. Later, L\'{a}sl\'{o} proved BW inequality for $3 \times 3$ real matrices \cite{Laslo}, Lu \cite{Lu1} and, independently, Wong and Jin \cite{Jin} proved it for $n \times n$ real matrices, after which the complex case was proved by B\"{o}ttcher and Wenzel in \cite{BW2008}. 

On the other hand, motivated by a problem in the field of open quantum systems \cite{Alicki}, we have introduced the following real-valued functional in \cite{ref:DFKO}: 
\begin{equation}\label{Rfunc}
r(A,B) = \frac{1}{2}\left( \<[B,A],BA\> +\<[B,A^*], BA^*\> \right)	
\end{equation}
for any $A,B \in M_n(\CA)$ where $A^*$ is the Hermitian conjugation of $A$ and $\< A,B\> = \tr(A^*B)$ is the Hilbert-Schmidt inner product. The sharp bound for the $r$-function was shown in \cite{ref:DFKO}:
\begin{equation}\label{eq:bdMain1}
	r(A,B) \le \frac{1+\sqrt{2}}{2} \|A\|^2 \|B\|^2. 
\end{equation}
This has been used to derive a universal constraint \cite{CKKS,K,KAW} for relaxation rates which is satisfied by any quantum Markovian dynamics described by completely positive dynamical semigroup, or GKLS master equation \cite{GKS}. 
For this reason, we called the functional \eqref{Rfunc} the $r$-function, where ``$r$'' stands for relaxation. 

In \cite{ref:DFKO}, we also observed a close relationship between BW inequality and inequality \eqref{eq:bdMain1}. If $B$ is normal, $r$-function reduces to the commutator:
\[
r(A,B) = \frac{1}{2} \|[A,B]\|^2.
\] 
Hence, the direct application of BW inequality leads the bound 
\[
r(A,B) \le \|A\|^2\|B\|^2. 
\] 
However, this is only true for a normal matrix $B$ and does not hold for general matrices. On the other hand, applying the triangle inequality, Schwarz inequality, the norm inequality $\|AB\| \le \|A\|\|B\|$ and BW inequality for \eqref{Rfunc} implies the following general bound:
\[
r(A,B) \le \sqrt{2} \|A\|^2\|B\|^2.  
\] 
However, this in turn is a loose inequality. In the end, the sharp bound for general matrices is given by \eqref{eq:bdMain1} and, therefore the relation between BW inequality and \eqref{eq:bdMain1} is not yet fully understood.  

The purpose of this paper is to provide one answer to this problem by introducing the following function (referred in this paper $f$-function) for complex matrices $A,B \in M_n(\CA)$ and one parameter $q \in \R$:
\begin{align}
f(A,B;q) &:= \Re \<[B,A] , [B,A]_q \> \nonumber \\
 &= \tr \left(A^\ct B^\ct B A + q B^\ct A^\ct A B - \frac{1+q}{2} (A^\ct B^\ct AB + B^\ct A^\ct B A)\right), \label{eq:f}
\end{align}
where $[A,B]_q:= AB - q BA$ is the $q$-deformed commutator \cite{J, M}.
When $q=1$, the $f$-function is exactly the norm squared of the commutator: 
\begin{equation}\label{eq:fNS}
f(A,B;1)= \|[A,B]\|^2.
\end{equation}
(See also \cite{CKOS} for another generalization $\|[A,B]_q\|^2$ and its norm bound.)  
When $q=0$, the $f$-function and $r$-function have the following relation:
\begin{equation}\label{eq:rf}
r(A,B) = \frac{f(A,B;0) + f(A^\ct,B;0)}{2}. 
\end{equation}
In this sense, the $f$-function connects the $r$-function and the norm squared of the commutator through one parameter $q$.

We consider the sharp bound for the $f$-function, and obtain a partial answer.

\begin{figure}[h]
\includegraphics[width=12cm]{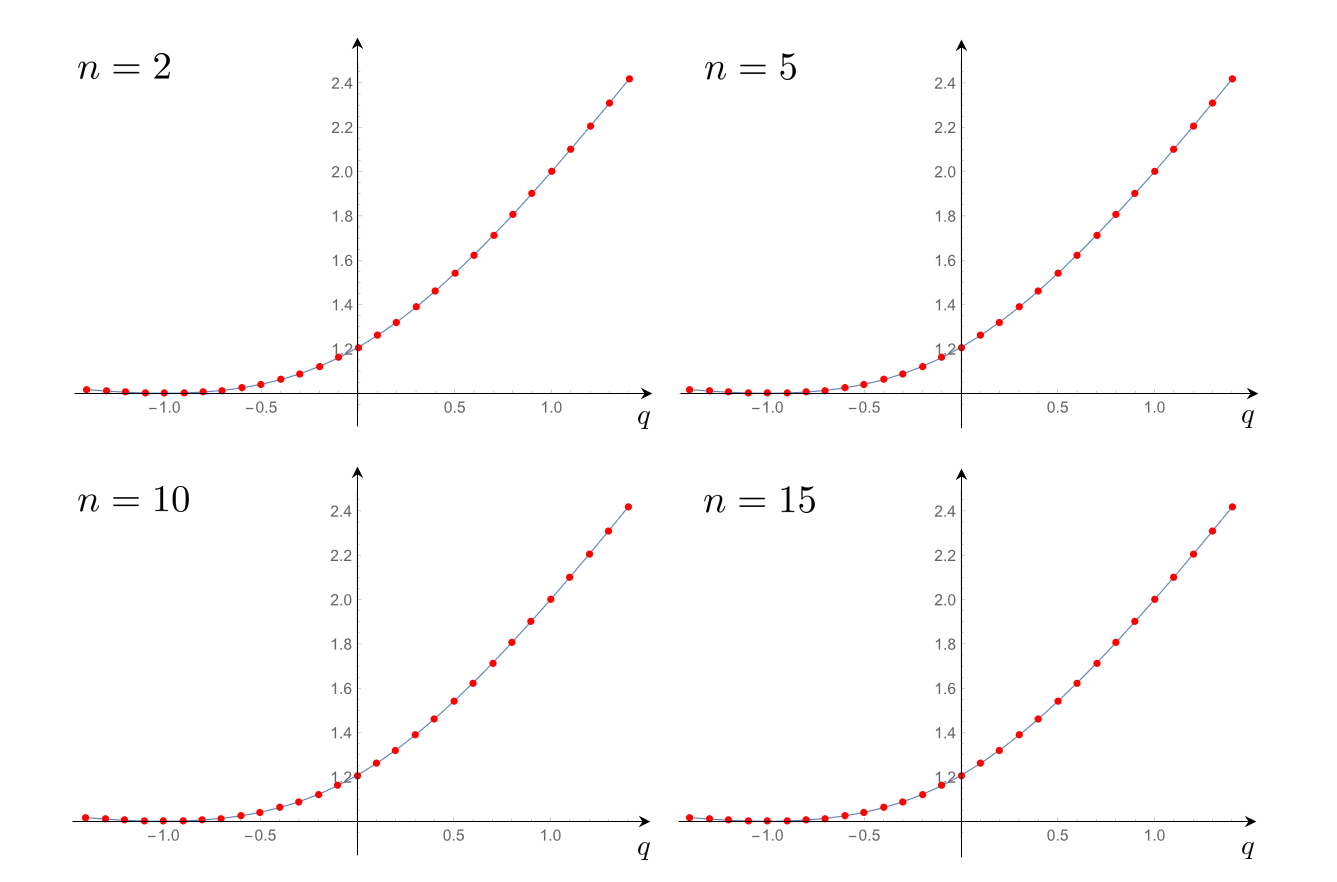}
\caption{Numerical evidences of \eqref{eq:Conj} for $n=2,5,10,15$: 
By the quadratic nature of the $f$-function, our conjecture \eqref{eq:Conj} can be rephrased as $\sup_{A,B \neq 0 \in M_n(\CA) } \frac{f(A,B;q)}{\|A|^2\|B\|^2} = c(q)$. 
Numerically optimized values (Red points) of the left hand side perfectly coincide with $c(q) =  ((1+q) + \sqrt{2(1+q^2)})/2$ (blue curves). }\label{fig:no}
\end{figure}

\begin{Thm}\label{thm:main}
Let $A$ and $B$ be $n\times n$ complex matrices.
If $A$ is a normal matrix or $n=2$, then the inequality
\begin{equation}\label{eq:Conj}
f(A,B;q) \le c(q) \|A\|^2 \|B\|^2
\end{equation}
holds, where 
\[
c(q) = \frac{(1+q) + \sqrt{2(1+q^2)}}{2}. 
\]
Moreover, this inequality is sharp, i.e., there are non-zero matrices $A$ and $B$ at which the bound is attained.  
\end{Thm}
We conjecture that the theorem is also true for any matrices $A, B \in M_n(\CA)$. This is numerically supported for $n$ up to $15$ (See Fig.~\ref{fig:no}).

Importantly, the inequality \eqref{eq:Conj} integrates BW inequality and the inequality \eqref{eq:bdMain1} as follows. 
First, inequality \eqref{eq:Conj} with $q = 1$ recovers BW inequality by noting $c(1) = 2$ and using relation \eqref{eq:fNS}.  
Hence, inequality \eqref{eq:Conj} is considered as a one-parameter generalization of BW inequality. 
Second, when $q=0$, we have $c(0) = \frac{1+\sqrt{2}}{2}$, and therefore, \eqref{eq:Conj} and \eqref{eq:rf} imply inequality \eqref{eq:bdMain1} for $r$-function. 
We also note that, the conjecture \eqref{eq:Conj} for $q=-1$ is easily shown by noting $c(-1) = 1$ and $f(A,B;-1) = \|B A\|^2 - \|A B\|^2$ and using the norm inequality $\|BA\| \le \|B\|\|A\|$.

\section{Proofs of the main theorem}


First, we mention some general facts for the $f$-function and $c(q)$.

\begin{Prop}
For any unitary matrix $U \in M_n(\CA)$ and complex numbers $\alpha, \beta \in \CA$,
the following equations hold:
\begin{align*}
f(UAU^\ct,UBU^\ct;q) &= f(A,B;q), \\
f(\alpha A,\beta B;q) &= |\alpha|^2 |\beta|^2 f(A,B;q).
\end{align*}
\end{Prop}
The proof is obvious and will be omitted.
This proposition says that the $f$-function is unitarily invariant.

In the following, we will write $c(q)=c$ unless there is no confusion.
Note that $c$ is bigger than or equal to $1$, $q$ and $1+q$ for any $q \in \R$.

\begin{Lem}\label{lem:lem1}
Let $\varepsilon_1= \varepsilon_1(q)$ and $\varepsilon_2= \varepsilon_2(q)$ be
\[
\varepsilon_1=
\begin{cases}
1 & q \ge -1 \\
-1 & q <-1
\end{cases}\quad {\rm and} \quad 
\varepsilon_2=
\begin{cases}
1 & q \le 1 \\
-1 & q > 1
\end{cases}.
\]
Then, the equations
\[
2\varepsilon_1 \sqrt{(c-q)(c-1)}  = (q+1) \quad {\rm and} \quad
2\varepsilon_2 \sqrt{ c(c-q-1)} = (1-q)
\]
hold.
\end{Lem}
[Proof] Since $\varepsilon_i^2 =1$ $(i=1,2)$, the solution for these equations is $c= \frac{(1+q) + \sqrt{2(1+q^2)}}{2}$.
\hfill $\blacksquare$

\medskip

Let $A = (a_{ij})_{i,j=1}^n$ and $B=(b_{ij})_{i,j=1}^n$ be arbitrary $n \times n$ matrices. 
A straightforward calculation shows
\begin{align*}
& c \|A\|^2\|B\|^2 -f(A,B;q) \\
& =
c\|A\|^2\|B\|^2 -\sum_{i,j,k,\ell=1}^n \left( a_{ij} \ov{a_{kj}} b_{\ell i} \ov{b_{\ell k}}
+q a_{ij} \ov{a_{i\ell}}b_{jk}\ov{b_{\ell k}} - (1+q)\Re(a_{ij} \ov{a_{\ell k}}b_{jk}\ov{b_{i\ell}}) \right).
\end{align*}
We define sets of indices as
\begin{align*}
&D_0 =\{ (i,i,i,i) \colon 1\le i \le n\} ,\\
&D_1 = \{(i,j,j,i) \colon 1\le i,j \le n, i\neq j \},\\
&D_2 = \{ (i,j,k,i) \colon 1\le i,j,k \le n, j\neq k\},\\ 
&D_3 = \{(i,j,j,k) \colon 1\le i,j,k \le n, i\neq k\}  \quad {\rm and}  \\
&D_4 =\{(i,j,k,\ell) \colon 1\le i,j,k,\ell \le n \} \backslash (D_1 \cup D_2 \cup D_3) .
\end{align*}
Note that $D_1, D_2, D_3$ and $D_4$ are mutually disjoint and $D_0 \subset D_4$.
The term $a_{ij} \ov{a_{kj}} b_{\ell i} \ov{b_{\ell k}}$ has the form 
$|a_{**}|^2 |b_{**}|^2$, when $i=k$. In this case, 
\[
a_{ij} \ov{a_{ij}} b_{\ell i} \ov{b_{\ell i}}= |a_{ij}|^2 |b_{\ell i}|^2.
\]
Therefore,
\begin{align*}
\sum_{i,j,\ell}a_{ij} \ov{a_{ij}} b_{\ell i} \ov{b_{\ell i}}= \sum_{i,j,\ell} |a_{ij}|^2 |b_{\ell i}|^2
=\sum_{(i,j,k,\ell) \in D_0 \cup D_1 \cup D_2} |a_{ij}|^2|b_{k\ell}|^2.
\end{align*}
Similarly, the term $a_{ij} \ov{a_{i\ell}}b_{jk}\ov{b_{\ell k}}$ has the form $|a_{**}|^2 |b_{**}|^2$, when $j=\ell$. Then, 
\[
\sum_{i,j,k} a_{ij} \ov{a_{ij}}b_{jk}\ov{b_{j k}} = \sum_{i,j,k} |a_{ij}|^2 |b_{jk}|^2
=\sum_{(i,j,k,\ell) \in D_0 \cup D_1 \cup D_3} |a_{ij}|^2|b_{k\ell}|^2.
\]
The term $a_{ij} \ov{a_{\ell k}}b_{jk}\ov{b_{i\ell}}$ has the form $|a_{**}|^2 |b_{**}|^2$, only when $i= j=k=\ell$. 

\noindent
Combining these, we have
\begin{align*}
& c \|A\|^2\|B\|^2 -f(A,B;q) \\
& =
c\sum_{(i,j,k,\ell) \in D_4}  |a_{ij}|^2|b_{k\ell}|^2 +
(c-1)\sum_{(i,j,k,\ell) \in D_2}|a_{ij}|^2|b_{k\ell}|^2 \\
&+ (c-q)\sum_{(i,j,k,\ell) \in D_3} |a_{ij}|^2|b_{k\ell}|^2 +
(c-1-q)\sum_{(i,j,k,\ell) \in D_1} |a_{ij}|^2|b_{k\ell}|^2 \\
&-\sum_{\substack{i,j,k,\ell\\ i\neq k}} a_{ij} \ov{a_{kj}} b_{\ell i} \ov{b_{\ell k}}
-q\sum_{\substack{i,j,k,\ell\\ j\neq \ell}}  a_{ij} \ov{a_{i\ell}}b_{jk}\ov{b_{\ell k}} 
+(1+q)\sum_{\substack{i,j,k,\ell \\ {\rm except} \\ i=j=k=\ell}}
\Re(a_{ij} \ov{a_{\ell k}}b_{jk}\ov{b_{i\ell}}) .
\end{align*}

\begin{Lem}
The next equation holds:
\begin{align*}
&\sum_{\substack{i,k \\ i\neq k}}
|\sqrt{c-q} \sum_{j=1}^n a_{ij}b_{jk} + 
\varepsilon_1 \sqrt{c-1} \sum_{\ell=1}^n a_{\ell k}b_{i\ell} |^2  \\
&=
(c-1) \sum_{D_2}|a_{ij}|^2|b_{k\ell}|^2 +
(c-q)\sum_{D_3} |a_{ij}|^2|b_{k\ell}|^2 \\
&+
(c-1) \sum_{\substack{i,j,k,l \\ i\neq k, j\neq \ell}} a_{ij}\ov{a_{kj}}b_{\ell i}\ov{b_{\ell k}}
+(c-q) \sum_{\substack{i,j,k,l \\ i\neq k, j\neq \ell}} a_{ij}\ov{a_{i\ell}}b_{jk}\ov{b_{\ell k}}\\
&+  (1+q)\Re\sum_{\substack{i,j,k,l \\ i\neq k}} a_{ij}\ov{a_{\ell k}}b_{jk}\ov{b_{i\ell}}.
\end{align*}
\end{Lem}
[Proof]  
Using the equation $\varepsilon_1 2\sqrt{(c-q)(c-1)} =(1+q)$ obtained in Lemma \ref{lem:lem1}, we can prove this by straightforward calculation.
\hfill $\blacksquare$ 

By this lemma, we can calculate $c \|A\|^2\|B\|^2 -f(A,B;q)$ as
\begin{align}
& c \|A\|^2\|B\|^2 -f(A,B;q)  \nonumber \\
& =
\sum_{\substack{i,k \\ i\neq k}}
|\sqrt{c-q} \sum_{j} a_{ij}b_{jk} +
\varepsilon_1 \sqrt{c-1} \sum_{\ell} a_{\ell k}b_{i\ell} |^2  \nonumber \\
&+c\sum_{D_4}  |a_{ij}|^2|b_{k\ell}|^2 +
(c-1-q)\sum_{ D_1} |a_{ij}|^2|b_{k\ell}|^2  \nonumber \\
& -c \sum_{\substack{i,j,k,\ell \\ i\neq k, j\neq \ell}} (a_{ij}\ov{a_{kj}}b_{\ell i}\ov{b_{\ell k}}
+a_{ij}\ov{a_{i\ell}}b_{jk}\ov{b_{\ell k}}) \nonumber  \\
&-\sum_{\substack{i,j,k \\ i\neq k}} a_{ij} \ov{a_{kj}} b_{j i} \ov{b_{j k}}
-q\sum_{\substack{i,j,\ell\\ j\neq \ell}}  a_{ij} \ov{a_{i\ell}}b_{ji}\ov{b_{\ell i}} 
+(1+q) \sum_{\substack{i,j,\ell\\ {\rm except} \\ i=j=\ell}} 
\Re(a_{ij} \ov{a_{\ell i}}b_{ji}\ov{b_{i\ell}}). \nonumber
\end{align}

\subsection{Proof for $n=2$}

In this subsection, the proof of Theorem \ref{thm:main} for $n=2$ is given. When $n=2$, a direct calculation shows: 
\begin{align*}
& | a_{11}\ov{b_{22}} + a_{22}\ov{b_{11}} 
-a_{12}\ov{b_{12}}-a_{21}\ov{b_{21}}|^2 \\
&=
\sum_{D_4\backslash D_0} |a_{ij}|^2 |b_{k\ell}|^2
-\sum_{\substack{i,j,k,\ell \\ i\neq k, j\neq \ell}} 
\left(a_{ij}\ov{a_{kj}}b_{\ell i}\ov{b_{\ell k}}
+a_{ij}\ov{a_{i\ell}}b_{jk}\ov{b_{\ell k}}
-a_{ij}\ov{a_{k\ell}}b_{ji}\ov{b_{\ell k}} \right).
\end{align*}
Using this, we have  
\begin{align}
& c \|A\|^2\|B\|^2 -f(A,B;q) \nonumber \\
& =
\sum_{\substack{i,k \\ i\neq k}}
|\sqrt{c-q} \sum_{j} a_{ij}b_{jk} +
\varepsilon_1 \sqrt{c-1} \sum_{\ell} a_{\ell k}b_{i\ell} |^2 \nonumber \\
&+c | a_{11}\ov{b_{22}} + a_{22}\ov{b_{11}} -a_{12}\ov{b_{12}}-a_{21}\ov{b_{21}}|^2  \nonumber \\
&+c\sum_{D_0}  |a_{ij}|^2|b_{k\ell}|^2 +
(c-1-q)\sum_{ D_1} |a_{ij}|^2|b_{k\ell}|^2 
-c \sum_{\substack{i,j,k,\ell \\ i\neq k, j\neq \ell}} a_{ij}\ov{a_{k\ell}}b_{ji}\ov{b_{\ell k}}\nonumber  \\
&-\sum_{\substack{i,j,k \\ i\neq k}} a_{ij} \ov{a_{kj}} b_{j i} \ov{b_{j k}}
-q\sum_{\substack{i,j,\ell\\ j\neq \ell}}  a_{ij} \ov{a_{i\ell}}b_{ji}\ov{b_{\ell i}} 
+(1+q) \sum_{\substack{i,j,\ell\\ {\rm except} \\ i=j=\ell}} 
\Re(a_{ij} \ov{a_{\ell i}}b_{ji}\ov{b_{i\ell}}) . \label{eq:form2}
\end{align}
By \eqref{eq:form2} and 
$2\varepsilon_2\sqrt{c(c-q-1)} =1-q$ (see Lemma \ref{lem:lem1}), 
we obtain
\begin{align*}
&c \|A\|^2\|B\|^2 -f(A,B;q) 
-
\sum_{\substack{i,k \\ i\neq k}}
|\sqrt{c-q} \sum_{j} a_{ij}b_{jk} +
\varepsilon_1\sqrt{c-1} \sum_{\ell} a_{\ell k}b_{i\ell} |^2 \\
&- c | a_{11}\ov{b_{22}} + a_{22}\ov{b_{11}} -a_{12}\ov{b_{12}}-a_{21}\ov{b_{21}}|^2  \\
&= c\left( |a_{11}|^2 |b_{11}|^2 + |a_{22}|^2|b_{22}|^2  \right)
+(c-1-q) \left( |a_{12}|^2|b_{21}|^2 + |a_{21}^2|b_{12}|^2 \right) \\
& -2c \Re \left(a_{11}\ov{a_{22}}b_{11}\ov{b_{22}} \right)  
-2(c-1-q) \Re \left(a_{12}\ov{a_{21}}b_{21}\ov{b_{12}} \right)   \\
&+(1-q) \Re \left( -a_{11}\ov{a_{21}}b_{11}\ov{b_{12}} -a_{12}\ov{a_{22}}b_{21}\ov{b_{22}}  
+a_{11}\ov{a_{12}}b_{11}\ov{b_{21}} +a_{21}\ov{a_{22}}b_{12}\ov{b_{22}} \right) \\
&=
\left|\sqrt{c}a_{11}b_{11} -\sqrt{c}a_{22}b_{22} 
+ 
\varepsilon_2 \left(\sqrt{c-1-q}a_{12}b_{21}-\sqrt{c-1-q}a_{21}b_{12} \right) \right|^2
\end{align*}
which implies the inequality \eqref{eq:Conj}.

To see that the inequality is sharp, one can easily find non-zero matrices $A$ and $B$ so that all the inequalities above are attained. 
For instance, let $a_{11} = \sqrt{c-1}$, $a_{22} = - \varepsilon_1\sqrt{c-q}$, $b_{12} = 1$ and the others are zero.
\hfill $\blacksquare$ 

\subsection{Proof for normal matrix}

In this subsection, the proof of Theorem \ref{thm:main} for normal matrix is given. Suppose that $A$ is a normal matrix. By the unitary invariance of the $f$-function, we can assume that $A$ is a diagonal matrix. Letting $A = {\rm diag}[a_1,a_2,\ldots,a_n]$ and $B = (b_{ij})$, direct computation shows  $\|A\|^2 = \sum_i |a_i|^2, \|B\|^2 = \sum_{i,j} |b_{ji}|^2$ and
\[
f(A,B;q) = \sum_{i \neq j}|b_{ji}|^2 (|a_i|^2 + q |a_j|^2 - (1+q) \Re(a_i \ov{a_j})).
\]
Hence, we have 
\begin{align*}
&c\|A\|^2\|B\|^2 - f(A,B;q) \\
&\ge \sum_{i\neq j} |b_{ji}|^2 (c (\sum_k |a_k|^2) - |a_i|^2 - q |a_j|^2 + (1+q) \Re(a_i \ov{a_j}) ) \\
&\ge \sum_{i\neq j} |b_{ji}|^2 ((c-1) |a_i|^2) +(c- q) |a_j|^2 + (1+q) \Re(a_i \ov{a_j}) ) \\
&= \sum_{i\neq j} |b_{ji}|^2 |\sqrt{c-1} a_i + 
\varepsilon_1 \sqrt{c- q} a_j |^2 \ge 0, 
\end{align*}
where in the equality we have used Lemma \ref{lem:lem1}. 
Therefore, the bound \eqref{eq:Conj} is shown for a normal matrix $A$ where $B$ is any complex matrix. 

One can prove that the inequality is sharp using the same example in the proof for $n=2$.
\hfill $\blacksquare$

\section{Conclusion and discussion}

In this paper, we have introduced the $f$-function \eqref{eq:f} which is a one parameter generalization of the norm square of the commutator and also $r$-function. The conjectured bound \eqref{eq:Conj} for the $f$-function integrates BW inequality and our previously obtained inequality \eqref{eq:bdMain1}. The conjecture is numerically supported for $n$ up to $15$, and is proved when one of the matrix is normal, or for general $2 \times 2$ matrices.

\bigskip

\noindent {\bf Acknowledgements}

\bigskip 

D.C. was supported by the Polish National Science Centre Project No. 2018/30/A/ST2/00837.


\begin{thebibliography}{plain}
\bibitem{Alicki} R. Alicki, K. Lendi, Quantum Dynamical Semigroups and Applications, Springer, Berlin, 1987.

\bibitem{ref:BW} A. B\"{o}ttcher, D. Wenzel, How big can the commutator of two matrices be and how big is it typically?, Linear Algebra Appl. {\bf 403}, (2005) 216-228.

\bibitem{BW2008} A. B\"{o}ttcher, D. Wenzel, The Frobenius norm and the commutator, Linear Algebra Appl. {\bf 429}, (2008) 1864-1185.


\bibitem{ref:DFKO}
D. Chru\'sci\'nski, R. Fujii, G. Kimura, H. Ohno, Constraints for the spectra of generators of quantum dynamical semigroups,  Linear Algebra Appl. {\bf 630}, (2021) 293-305. 

\bibitem{CKOS}
D. Chru\'sci\'nski, G. Kimura, H. Ohno, T. Singal, Bounding the Frobenius norm of a $q$-deformed commutator, Linear Algebra Appl. {\bf 646}, (2022) 95-106.

\bibitem{CKKS} D. Chruscinski, G. Kimura, A. Kossakowski, Y. Shishido, On the universal constraints for relaxation rates for quantum dynamical semigroup, Phys. Rev. Lett. {\bf 127},  (2021) 050401. 

\bibitem{GKS} V. Gorini, A. Kossakowski, E. C. G. Sudarshan, Completely positive dynamical semigroups of $N$-level systems, J. Math. Phys. {\bf 17}, (1976) 821. 

\bibitem{J} 
R. Jagannathan, Some introductory notes on quantum groups, quantum algebras, and their appli-cations, arXiv:math-ph/0105002.

\bibitem{ref:L}  G. Lindblad, On the generators of quantum dynamical semigroups, Comm. Math. Phys. {\bf 48},  (1976) 119-130.

\bibitem{K} G. Kimura, Restriction on relaxation times derived from the Lindblad-type master equations for two-level systems, Phys. Rev. A {\bf 66}, (2002) 062113.

\bibitem{KAW} G. Kimura, S. Ajisaka, K. Watanabe, Universal Constraints on Relaxation Times for $d$-Level GKLS Master Equations, Open Syst. Inform. Dynam. {\bf 24}, 
(2017) 1740009.

\bibitem{Laslo} L. Laslo, Proof of B\"{o}ttcher and Wenzel’s conjecture on commutator norms for 3-by-3 matrices, Linear Algebra Appl. {\bf 422}, (2007) 659–663.

\bibitem{Lu1}  Z. Lu, Normal Scalar Curvature Conjecture and its applications, Journal of Functional Analysis {\bf 261}, (2011) 1284–1308.

\bibitem{M}
S. Majid, Foundations of Quantum Group Theory, Cambridge University Press, Cambridge, 1995.

\bibitem{Jin}  S.~W. Vong, X.~Q. Jin, ``Proof of B\"{o}ttcher and Wenzel's conjecture", Oper. Matrices {\bf 2}, (2008) 435–442.

\end{thebibliography}
\end{document}